\documentclass[final,3p,times]{elsarticle}
\usepackage[latin9]{inputenc}
\usepackage{bm}
\usepackage{amsmath}
\usepackage{amssymb}
\usepackage{graphicx}
\usepackage{esint}
\usepackage{multirow}
\usepackage[subfigure]{graphfig}
\usepackage{rotating}
\usepackage{verbatim}
\usepackage{setspace}
\makeatletter
\newcommand*{\rom}[1]{\expandafter\@slowromancap\romannumeral #1@}

\usepackage{tabularx}\usepackage{subfigure}\usepackage{subfigure}\usepackage{color}

\usepackage{setspace}

\usepackage{bm}

\journal{}

\makeatother

\color{black}

\begin{document}
\title{Some practical versions of boundary variation diminishing (BVD) algorithm}

\author[ad1]{Xi Deng}

\author[ad1]{Bin Xie}

\author[ad1]{Feng Xiao \corref{cor}}


\address[ad1]{Department of Mechanical Engineering, Tokyo Institute of Technology, \\
	4259 Nagatsuta Midori-ku, Yokohama, 226-8502, Japan.}
\cortext[cor]{Corresponding author: Dr. F. Xiao (Email: xiao@es.titech.ac.jp)}

\begin{abstract}
This short note presents some variant schemes of boundary variation diminishing (BVD) algorithm in one dimension with the results of numerical tests for linear advection equation to facilitate practical use.  In spite of being presented in 1D fashion, all the schemes are simple and easy to  implement in multi-dimensions on structured and unstructured grids for nonlinear and system equations. 
\end{abstract}

\maketitle

\section{Introduction}
Spatial reconstruction schemes under the frameworks of finite volume method (FVM) or finite difference method (FDM) for hyperbolic conservation laws have been extensively exploited in the past decades using high-order-polynomial fashion interpolation functions, see \cite{abgrall-shu2016} for comprehensive and updated review on the developments in the fields. The polynomial-based high-order (PHO) reconstructions are highly demanded and show excellent performance in resolving numerical solutions which have relatively smooth variations in space, such as acoustic waves and vortices in compressible flow. However, when applied to strong discontinuities, the PHO reconstructions have to be projected to lower order or smoother polynomials (or more precisely rational functions in nonlinear schemes) to suppress numerical oscillations, which is called limiting projection. The limiting projection can be designed by following principles demanded from mathematical or physical perspectives, like the TVD, ENO/WENO and maximum-principle-satisfying. The PHO reconstructions with limiting projection might result in excessive numerical dissipation and tend to smear out the discontinuities in numerical solution. Being aware that the polynomials might not be the best choice in some circumstances for constructing the numerical solution, \cite{Marquina94, xiao04, Schroll06} suggested other non-polynomial functions with better monotonicity-preserving property, like  piecewise hyperbolic and  piecewise rational reconstructions, which over-perform the polynomial reconstructions in capturing discontinuities with monotone distributions and are able to provide oscillation-less solutions even without limiting projections. However, extending this sort of methods to higher order seems to be not straightforward.  

We have shown recently that combing the PHO reconstructions with Sigmoid functions provides a very promising methodology to construct high-fidelity schemes to capture both smooth and jump-like solutions \cite{Sun, Xie,deng2017-0,deng2017}. Implementing such an idea requires some principles or guidelines to hybridize (or switch) the candidate reconstructions according to the numerical solutions. Among other options, we have practiced the boundary variation diminishing (BVD) principle \cite{Sun} that  minimizes  the variations (jumps) of the reconstructed values at cell boundaries. The BVD reconstruction effectively reduces the dissipation in numerical fluxes, and thus results in numerical solutions with minimized numerical oscillation and dissipation errors. We have experimented the BVD principle with the 5th-order WENO scheme \cite{wenoz} and the THINC (Tangent of Hyperbola for INterface Capturing) scheme \cite{xiao_thinc2,xiao_thinc,ii2014,xie2014}, and found that the numerical results of high quality for both smooth and non-smooth solutions can be obtained. 

In order to facilitate the practical use of the BVD principle, as well as to promote further exploration for new schemes of better properties, in this note we summarize some simple and efficient BVD algorithms we have tested. Although the discussions are limited to 1D with numerical results only for linear advection test, all algorithms presented here as reconstruction schemes are applicable straightforwardly to multi-dimensions and nonlinear system equations.

\section{Preliminary \label{sec:model}}
In this paper, the 1D scalar conservation law in following form is used to introduce the BVD algorithm
\begin{equation}
\label{eq:scalar}
\frac{\partial q}{\partial t} + \frac{\partial f(q)}{\partial x} = 0,
\end{equation}
where $q(x,t)$ is the solution function and $f(q)$ is the flux function. We divide the computational domain into $N$ non-overlapping cell elements, ${\mathcal I}_{i}: x \in [x_{i-1/2},x_{i+1/2} ]$, $i=1,2,\ldots,N$, with a uniform grid spacing $\Delta x=x_{i+1/2}-x_{i-1/2}$. For a standard finite volume method, the volume-integrated average value $\bar{q}_{i}(t)$ in cell ${\mathcal I}_{i}$ is defined as
\begin{equation}
\bar{q}_{i}(t) 
\approx \frac{1}{\Delta x} \int_{x_{i-1/2}}^{x_{i+1/2}}
q(x,t) \; dx.
\end{equation}
The semi-discrete version of Eq.~(\ref{eq:scalar}) in the finite volume form can be expressed as a ordinary differential equations (ODEs)
\begin{equation}
\frac{\partial \bar{q}(t)}{\partial t}  =-\frac{1}{\Delta x}(\tilde{f}_{i+1/2}-\tilde{f}_{i-1/2}),
\end{equation}
where the numerical fluxes $\tilde{f}$ at cell boundaries can be computed by a Riemann solver
\begin{equation}
\tilde{f}_{i+1/2}=f_{i+1/2}^{\text{Riemann}}(q_{i+1/2}^{L},q_{i+1/2}^{R})
\end{equation}
if the reconstructed left-side value $q_{i+1/2}^{L}$ and right-side value $q_{i+1/2}^{R}$ at cell boundaries are provided. Essentially, the Riemann flux can be written in a canonical form as
\begin{equation}
\label{eq:Riemann}
f_{i+1/2}^{\text{Riemann}}(q_{i+1/2}^{L},q_{i+1/2}^{R}) =\frac{1}{2}\left(f(q_{i+1/2}^{L})+f(q_{i+1/2}^{R})\right)-\frac{|a_{i+1/2}|}{2}\left(q_{i+1/2}^{R}-q_{i+1/2}^{L})\right),
\end{equation}
where $a_{i+1/2}$ stands for the characteristic speed of the hyperbolic conservation law. The remaining main task is how to calculate $q_{i+1/2}^{L}$ and $q_{i+1/2}^{R}$ through the reconstruction process.

An effective and perhaps the most popular reconstruction scheme is the WENO schemes in which high order but non-oscillatory interpolation is achieved by a combination of several lower degree polynomials. Different variants have been devised following the pioneer works in \cite{liu,Jiang}. In present work, we use the improved  5th-order WENO scheme (WENOZ scheme) proposed in  \cite{wenoz} as the PHO reconstruction candidate in the BVD algorithm. Referring the interested readers to \cite{wenoz} for algorithmic details, we denote the 5th-order WENOZ reconstruction function  in cell ${\mathcal I}_{i}$ by ${q}^{W}_{i}(x)$, and the  cell boundary values are thus computed by  ${q}^{W}_{i}(x_{i-1/2})$ and ${q}^{W}_{i}(x_{i+1/2})$ respectively.

We use the THINC scheme \cite{xiao_thinc,xiao_thinc2} as another candidate for reconstruction. Being a Sigmod-type function, the THINC reconstruction (hyperbolic tangent function) is a differentiable and monotone function that fits well a step-like discontinuity. The THINC reconstruction function is written as
\begin{equation}
{q}^{T}_{i}(x)=\bar{q}_{min}+\dfrac{\bar{q}_{max}}{2} \left(1+\theta~\tanh \left(\beta \left(\dfrac{x-x_{i-1/2}}{x_{i+1/2}-x_{i-1/2}}-\tilde{x}_{i}\right)\right)\right),
\end{equation} 
where $\bar{q}_{min}=min(\bar{q}_{i-1},\bar{q}_{i+1})$, $\bar{q}_{max}=max(\bar{q}_{i-1},\bar{q}_{i+1})-\bar{q}_{min}$ and $\theta=sgn(\bar{q}_{i+1}-\bar{q}_{i-1})$. The jump thickness is controlled by parameter $\beta$. In our numerical tests shown later a constant value of $\beta=1.8$ is used, or explicitly stated otherwise. The unknown $\tilde{x}_{i}$, which represents the location of the jump center, is computed from $ \bar{q}_{i} = \frac{1}{\Delta x} \int_{x_{i-1/2}}^{x_{i+1/2}} {q}_{i}(x)^{T} \; dx$. Then the reconstructed values at cell boundaries by THINC function can be expressed by 
\begin{equation}
\begin{aligned}
&{q}^{T}_{i}(x_{i+1/2})=\bar{q}_{min}+\dfrac{\bar{q}_{max}}{2} \left(1+\theta \dfrac{\tanh(\beta)+A}{1+A~\tanh(\beta)}\right)\\
&{q}^{T}_{i}(x_{i-1/2})=\bar{q}_{min}+\dfrac{\bar{q}_{max}}{2} \left(1+\theta~ A\right)
\end{aligned}
\end{equation}
where $A=\frac{B/\cosh(\beta)-1}{\tanh(\beta)}$, $B=\exp(\theta~\beta(2~C-1))$ and $C=\dfrac{\bar{q}_{i}-\bar{q}_{min}+\epsilon}{\bar{q}_{max}+\epsilon}$  with $\epsilon=10^{-20}$.

\section{The BVD algorithms \label{bvd}}  

Given two reconstruction functions shown above, i.e. ${q}^{W}_{i}(x)$ and ${q}^{T}_{i}(x)$, we use the BVD principle to choose the final reconstruction so that  $|q_{i+1/2}^{R}-q_{i+1/2}^{L}|$ is minimized, which then effectively reduces the numerical dissipation, i.e. the 2nd  term in Eq.~(\ref{eq:Riemann}). In this note, we present several BVD algorithms that compare the reconstructed values across cell boundaries. It is noted that we are focusing on the simple and easy-to-use versions of BVD, which are called in turn BVD(\rom{1}), BVD(\rom{2}), BVD(\rom{3}) and BVD(\rom{4}) in this note. Using  WENO and THINC as the candidate reconstructions, the resulted schemes are called WENO-THINC-BVD(\rom{1}) $\sim$ (\rom{4}) methods.

\subsection{BVD(\rom{1})  algorithm\cite{Sun}}

\begin{enumerate}[i)]
\item Find $q^{\xi}_i(x)$ and $q^{\eta}_{i+1}(x)$ with $\xi$ and $\eta$ being either $W$ or $T$, so that the boundary variation (BV)
\begin{equation} \label{min_ip12}
BV (q)_{i+\frac{1}{2}} =|q^{\xi}_i(x_{i+\frac{1}{2}})-q^{\eta}_{i+1}(x_{i+\frac{1}{2}})|,
\end{equation}
is minimized;
\item In case that a different choice for cell ${\mathcal I}_{i}$ is made when applying step i) to the neighboring interface $x_{i-\frac{1}{2}}$, that is,  $q^{\xi'}_i(x)$ found to minimize 
\begin{equation}\label{min_im12}
BV (q)_{i-\frac{1}{2}} =|q^{\xi'}_{i-1}(x_{i-\frac{1}{2}})-q^{\eta'}_{i}(x_{i-\frac{1}{2}})|,
\end{equation}
with $\xi'$ and $\eta'$ being either $W$ or $T$, is different from $q^{\xi}_i(x)$ found to minimize \eqref{min_ip12}, adopt the following criterion to uniquely determine the reconstruction function. 
\begin{equation}
q_{i}(x)=
\left\{ \begin{array}{l}
q^{W}_{i}(x),  \  {\rm if} \ \left(q^{\xi}_i(x_{i+\frac{1}{2}})-q^{\eta}_{i+1}(x_{i+\frac{1}{2}})\right)\left(q^{\xi'}_{i-1}(x_{i-\frac{1}{2}})-q^{\eta'}_{i}(x_{i-\frac{1}{2}}\right) < 0, \\ 
q^{T}_{i}(x),  \  {\rm otherwise. }
\end{array}\right.
\label{ifcase1}
\end{equation}

\item  Compute the left-side value $ u^{L}_{i+\frac{1}{2}} $ at $x_{i+\frac{1}{2}}$ and the right-side value $ u^{R}_{i-\frac{1}{2}} $ at $x_{i-\frac{1}{2}}$ by 
\begin{equation}
u^{L}_{i+\frac{1}{2}}=q_i(x_{i+\frac{1}{2}}) \ \ {\rm and} \ \  u^{R}_{i-\frac{1}{2}}=q_{i}(x_{i-\frac{1}{2}}).
\end{equation}

\end{enumerate}

\subsection{BVD(\rom{2})  algorithm\cite{deng2017}}

\begin{enumerate}[i)]
\item For each reconstruction  $q^{\xi}_i(x)$  with $\xi$ being either $W$ or $T$, calculate the minimum value of total boundary variation (TBV), $ mTBV^{\xi}_i$,  for all possible candidate reconstructions over the neighboring cells by 
 \begin{equation}
 \begin{aligned}
 mTBV^{\xi}_i=\min \bigg( &|q^{W}_{i-1}(x_{i-\frac{1}{2}})-q^{\xi}_i(x_{i-\frac{1}{2}})|+|q^{W}_{i+1}(x_{i+\frac{1}{2}})-q^{\xi}_i(x_{i+\frac{1}{2}})|,  \\
&|q^{T}_{i-1}(x_{i-\frac{1}{2}})-q^{\xi}_i(x_{i-\frac{1}{2}})|+|q^{T}_{i+1}(x_{i+\frac{1}{2}})-q^{\xi}_i(x_{i+\frac{1}{2}})|, \\
&|q^{W}_{i-1}(x_{i-\frac{1}{2}})-q^{\xi}_i(x_{i-\frac{1}{2}})|+|q^{T}_{i+1}(x_{i+\frac{1}{2}})-q^{\xi}_i(x_{i+\frac{1}{2}})|,  \\
&|q^{T}_{i-1}(x_{i-\frac{1}{2}})-q^{\xi}_i(x_{i-\frac{1}{2}})|+|q^{W}_{i+1}(x_{i+\frac{1}{2}})-q^{\xi}_i(x_{i+\frac{1}{2}})| \bigg).
 \end{aligned}
 \label{mintbv}
 \end{equation}  
  
 \item 
 Given the minimum TBVs for both  $q^{W}_i(x)$  and  $q^{T}_i(x)$, $ mTBV^{W}_i$ and $ mTBV^{T}_i$ computed from \eqref{mintbv},  
 choose the reconstruction function for cell ${\mathcal I}_{i}$ by 
 \begin{equation}
{q}_{i}(x)=\left\{
 \begin{array}{l}
{q}^{T}_{i},~~~{\rm if } \ \  mTBV^{T}_i <  mTBV^{W}_i, \\
{q}^{W}_{i}~~~~\mathrm{otherwise}
 \end{array}
 \right..
 \end{equation}
That is the THINC reconstruction function will be employed in the targeted cell $I_{i}$ if the minimum TBV value of THINC is smaller than that of high order WENO interpolation.

\item  Compute the left-side value $ u^{L}_{i+\frac{1}{2}} $ at $x_{i+\frac{1}{2}}$ and the right-side value $ u^{R}_{i-\frac{1}{2}} $ at $x_{i-\frac{1}{2}}$ by 
\begin{equation}
u^{L}_{i+\frac{1}{2}}=q_i(x_{i+\frac{1}{2}}) \ \ {\rm and} \ \  u^{R}_{i-\frac{1}{2}}=q_{i}(x_{i-\frac{1}{2}}).
\end{equation}

\end{enumerate}

\subsection{BVD(\rom{3})  algorithm \cite{Xie} } 

\begin{enumerate}[i)]
\item Compute the TBV of the target cell $I_{i}$ with the WENO reconstruction by
\begin{equation}
 TBV_{i}^{(W)}=\dfrac{\bigg(q^{W}_{i-1}(x_{i-\frac{1}{2}})-q^{W}_{i}(x_{i-\frac{1}{2}})\bigg)^4+\bigg(q^{W}_{i+1}(x_{i+\frac{1}{2}})-q^{W}_{i}(x_{i+\frac{1}{2}})\bigg)^4}{\big(\bar{q}_{i}-\bar{q}_{i-1}\big)^4+\big(\bar{q}_{i}-\bar{q}_{i+1}\big)^4+\epsilon}, 
\end{equation}   
where $\epsilon$ is a small positive of $10^{-16}$ for avoiding zero-division. 
\item Compute  the smoothness indicator by
\begin{equation}
S=\dfrac{1-TBV_{i}^{(W)}}{\max(TBV_{i}^{(W)},\epsilon)},
\end{equation}  
The cutoff number $S_{c}$ is used as the threshold value so that the cell where $S<S_{c}$ is identified to contain a non-smooth solution. We set $S_{c}=1\times10^{6}$ here. 

\item For cell $I_{i}$ where $S<S_{c}$, the reconstruction function is determined by blending ${q}_{i}(x)^{T}$ and ${q}_{i}(x)^{W}$ as follows 
\begin{equation}
{q}_{i}(x)=\omega_{i}{q}_{i}(x)^{T}+(1-\omega_{i}){q}_{i}(x)^{W},
\label{blending}
\end{equation}      
where $\omega_{i}$ is a weight parameter. By assuming that the WENO reconstruction is applied on neighbor cells, $\omega_{i}$ is obtained by minimizing
\begin{equation}
\tau_{i}=\bigg(q^{W}_{i-1}(x_{i-\frac{1}{2}})-q_{i}(x_{i-\frac{1}{2}})\bigg)^2+\bigg(q^{W}_{i+1}(x_{i+\frac{1}{2}})-q_{i}(x_{i+\frac{1}{2}})\bigg)^2, 
\end{equation} 
which leads to 
\begin{equation}
\dfrac{\partial \tau_{i}}{\partial \omega_{i}}=0.
\end{equation}
As long as $\omega_{i}$ is determined, the reconstruction function ${q}_{i}(x)$ is computed from \eqref{blending}. 

\item  Compute the left-side value $ u^{L}_{i+\frac{1}{2}} $ at $x_{i+\frac{1}{2}}$ and the right-side value $ u^{R}_{i-\frac{1}{2}} $ at $x_{i-\frac{1}{2}}$ by 
\begin{equation}
u^{L}_{i+\frac{1}{2}}=q_i(x_{i+\frac{1}{2}}) \ \ {\rm and} \ \  u^{R}_{i-\frac{1}{2}}=q_{i}(x_{i-\frac{1}{2}}).
\end{equation}

\end{enumerate}

This simplified variant can be easily used for unstructured grids as shown in \cite{Xie,deng2017-0, dengAIAA}. 

\subsection{BVD(\rom{4})  algorithm}

\begin{enumerate}[i)]
\item Compute the TBVs of the target cell $I_{i}$ using WENO and THINC for  $I_{i}$  and its two neighboring cells respectively, 
\begin{equation}
TBV_{i}^{W}=|q^{W}_{i-1}(x_{i-\frac{1}{2}})-q^{W}_{i}(x_{i-\frac{1}{2}})|+|q^{W}_{i}(x_{i+\frac{1}{2}})-q^{W}_{i+1}(x_{i+\frac{1}{2}})|
\end{equation}
and 
\begin{equation}
TBV_{i}^{T}=|q^{T}_{i-1}(x_{i-\frac{1}{2}})-q^{T}_{i}(x_{i-\frac{1}{2}})|+|q^{T}_{i}(x_{i+\frac{1}{2}})-q^{T}_{i+1}(x_{i+\frac{1}{2}})|. 
\end{equation}   

\item 
 Given TBVs for both  $q^{W}_i(x)$  and  $q^{T}_i(x)$, $TBV_{i}^{W}$ and $TBV_{i}^{T}$,  
 choose the reconstruction function for cell ${\mathcal I}_{i}$ by 
\begin{equation}
{q}_{i}(x)=\left\{
 \begin{array}{l}
{q}^{T}_{i}~~~{\rm if } \ \ TBV_{i}^{T}  < TBV_{i}^{W}, \\
{q}^{W}_{i}~~~~\mathrm{otherwise}
\end{array}
\right..
 \end{equation}
 
\item  Compute the left-side value $ u^{L}_{i+\frac{1}{2}} $ at $x_{i+\frac{1}{2}}$ and the right-side value $ u^{R}_{i-\frac{1}{2}} $ at $x_{i-\frac{1}{2}}$ by 
\begin{equation}
u^{L}_{i+\frac{1}{2}}=q_i(x_{i+\frac{1}{2}}) \ \ {\rm and} \ \  u^{R}_{i-\frac{1}{2}}=q_{i}(x_{i-\frac{1}{2}}).
\end{equation}

\end{enumerate}

\section{Numerical examples \label{sec:results}}

We show in this section the numerical results of above schemes for advection test. 
It is noted that when implementing the BVD algorithms presented above we make use the following criteria for choosing THINC reconstruction.   
\begin{equation}
\begin{aligned}
&\delta<C<1-\delta,\\
&(\bar{q}_{i+1}-\bar{q}_{i})(\bar{q}_{i}-\bar{q}_{i-1})>0,
\end{aligned}
\end{equation}
where $\delta$ is a small positive $10^{-4}$ in all numerical tests presented next. Otherwise, the WENO reconstruction is used.
 
\subsection{Advection of One-Dimensional Complex Waves}
\begin{figure}
	\subfigure{\centering\includegraphics[scale=0.5,trim={0.5cm 0.5cm 0.5cm 0.5cm},clip]{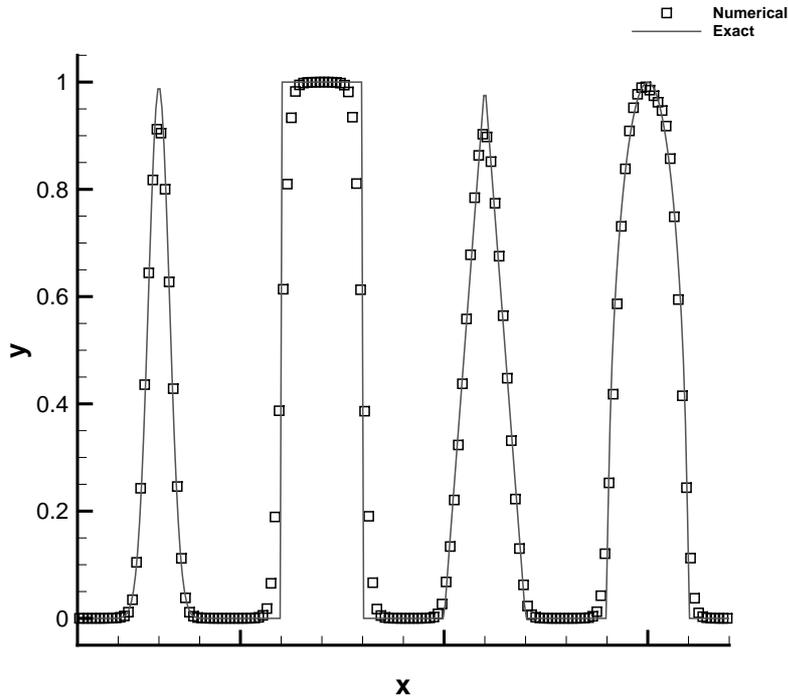}}
	\protect\caption{Numerical results of advection of complex waves with the WENOZ scheme.\label{fig:wenoz}}	
\end{figure}

Proposed in \cite{Jiang}, the test of propagation of a complex wave which includes both discontinuous and smooth solutions has been used widely to examine the performance of numerical schemes in solving profiles of different smoothness. The initial distribution of the advected field is set the same as \cite{Jiang}. The numerical result with the WENOZ scheme after one period of computation on a 200-cell mesh is plotted in Fig.~\ref{fig:wenoz}. Although WENOZ has good performance for smooth region, the discontinuity has been diffused by nearly 8 cells. The smeared discontinuity will become worse for a long time computation.  

\begin{figure}
	\subfigure{\centering\includegraphics[scale=0.5,trim={0.5cm 0.5cm 0.5cm 0.5cm},clip]{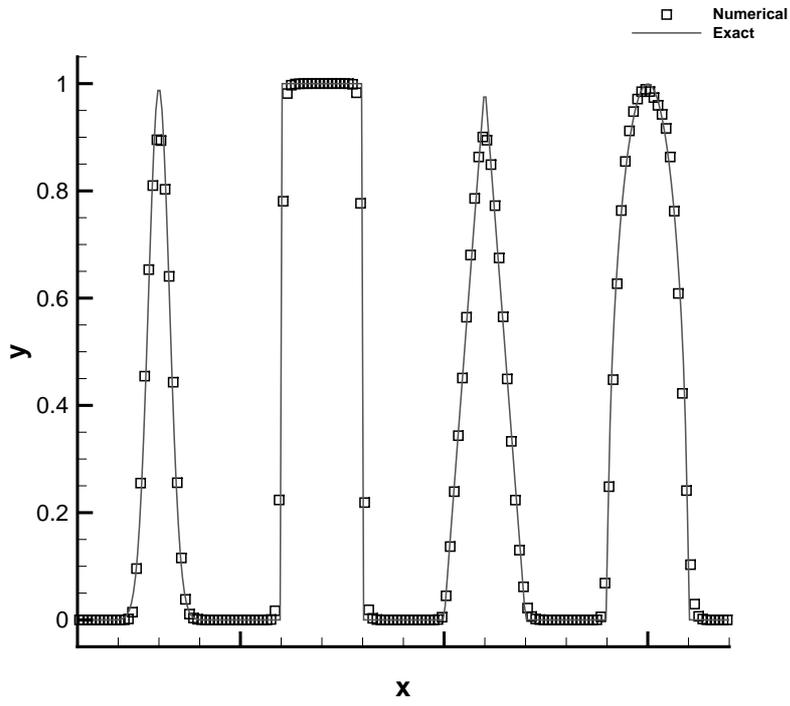}}
	\protect\caption{Numerical results of advection of complex waves with the WENOZ-THINC-BVD scheme with BVD(\rom{1}) algorithm.\label{fig:sunwenoz}}	
\end{figure}

\begin{figure}
	\subfigure{\centering\includegraphics[scale=0.5,trim={0.5cm 0.5cm 0.5cm 0.5cm},clip]{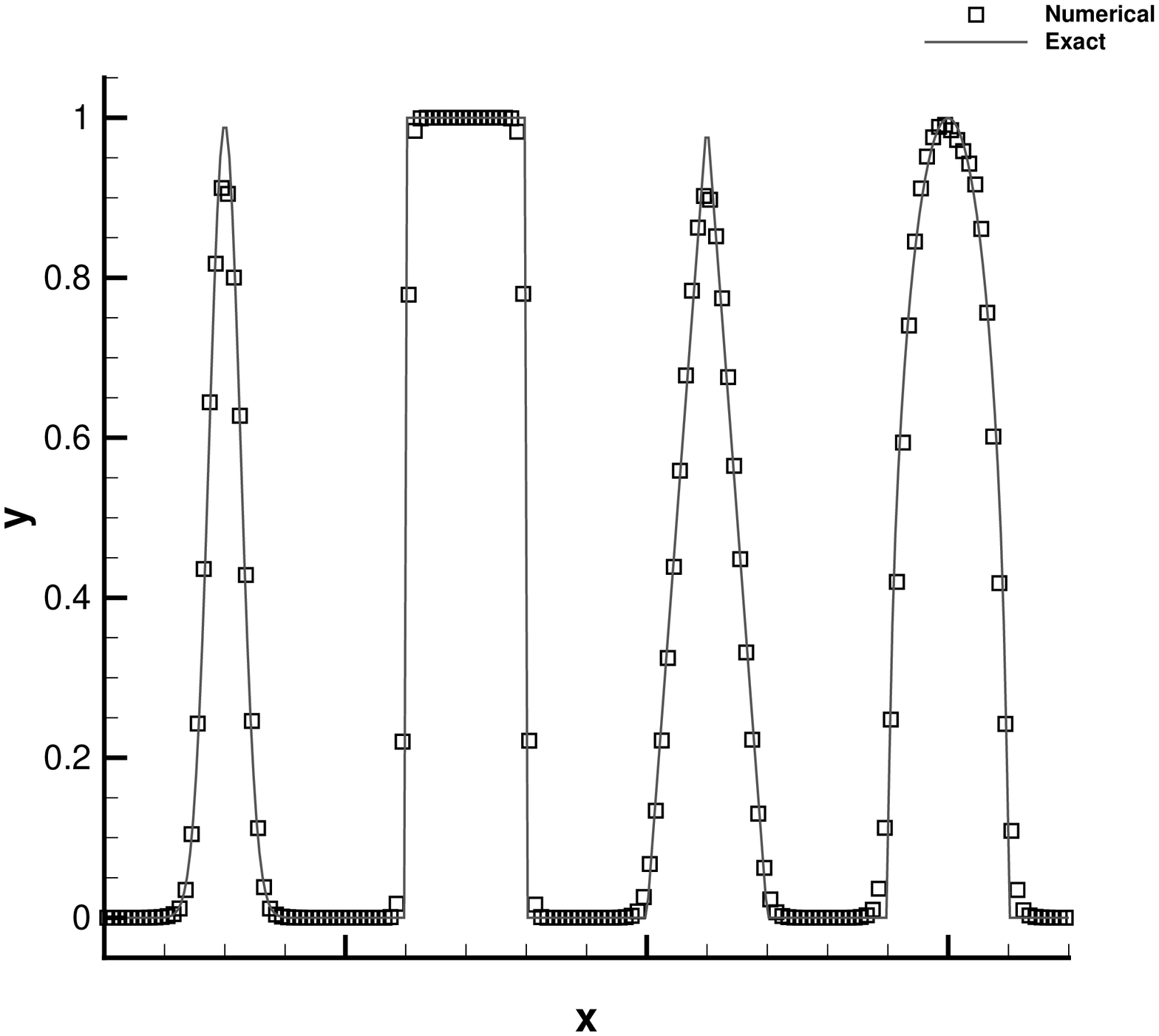}}
	\protect\caption{Same as Fig.\ref{fig:sunwenoz}, but with BVD(\rom{2}) algorithm.\label{fig:deng1wenoz}}	
\end{figure}

\begin{figure}
	\subfigure{\centering\includegraphics[scale=0.5,trim={0.5cm 0.5cm 0.5cm 0.5cm},clip]{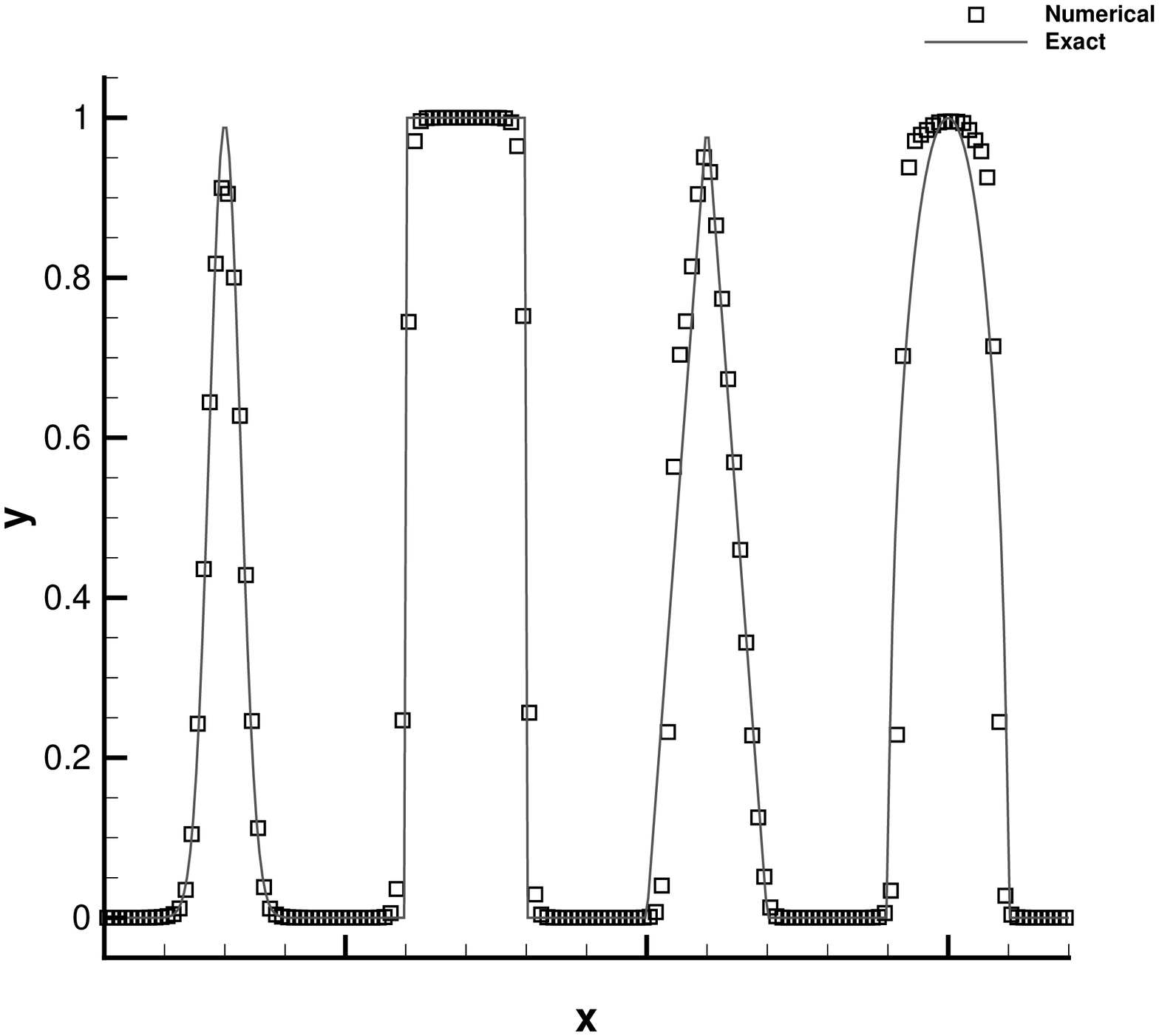}}
	\protect\caption{Same as Fig.\ref{fig:sunwenoz}, but with BVD(\rom{3}) algorithm.\label{fig:xiewenoz}}	
\end{figure}

\begin{figure}
	\subfigure{\centering\includegraphics[scale=0.5,trim={0.5cm 0.5cm 0.5cm 0.5cm},clip]{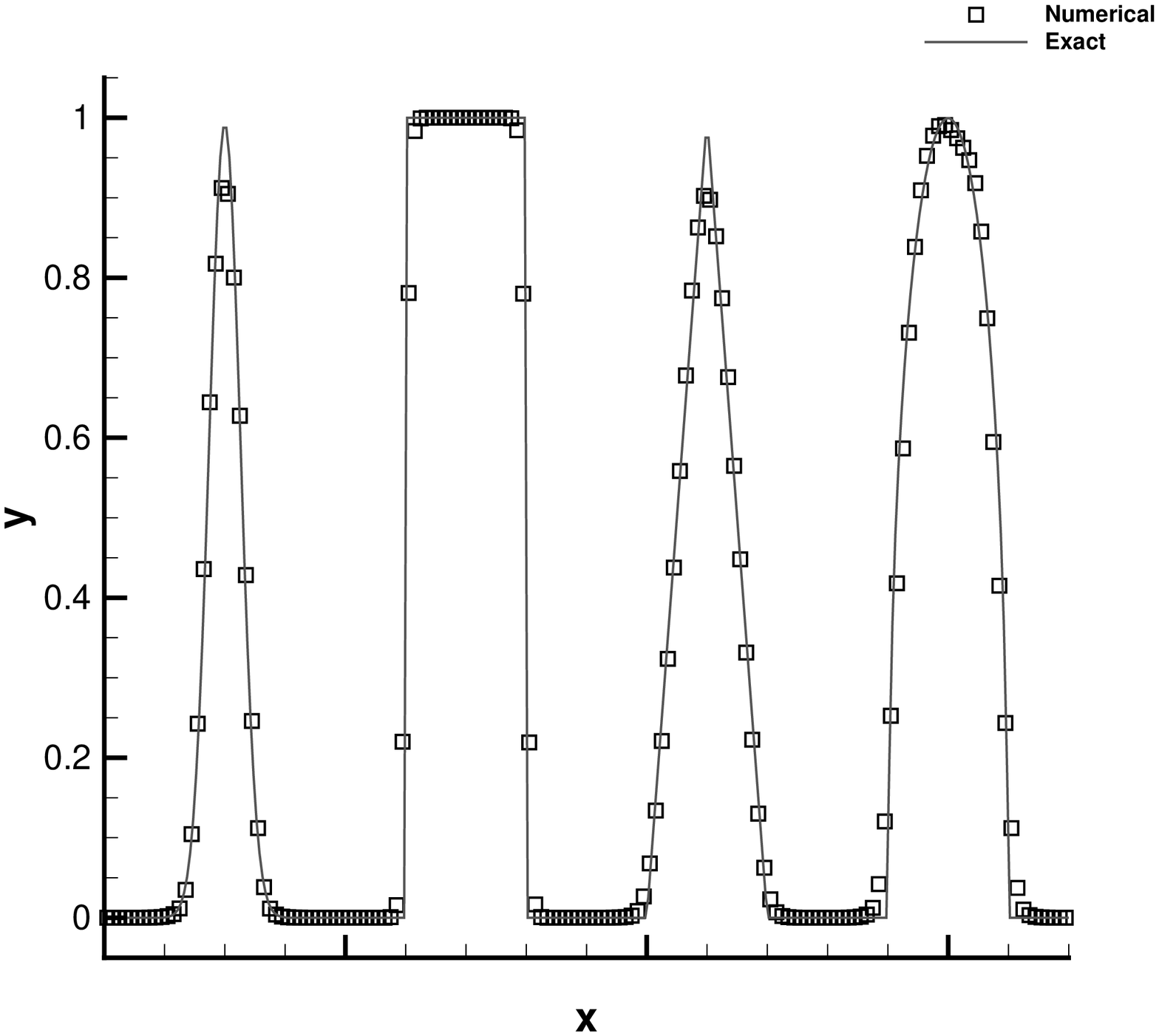}}
	\protect\caption{Same as Fig.\ref{fig:sunwenoz}, but with BVD(\rom{2}) algorithm.\label{fig:deng2wenoz}}	
\end{figure}

\begin{figure}
	\subfigure{\centering\includegraphics[scale=0.5,trim={0.5cm 0.5cm 0.5cm 0.5cm},clip]{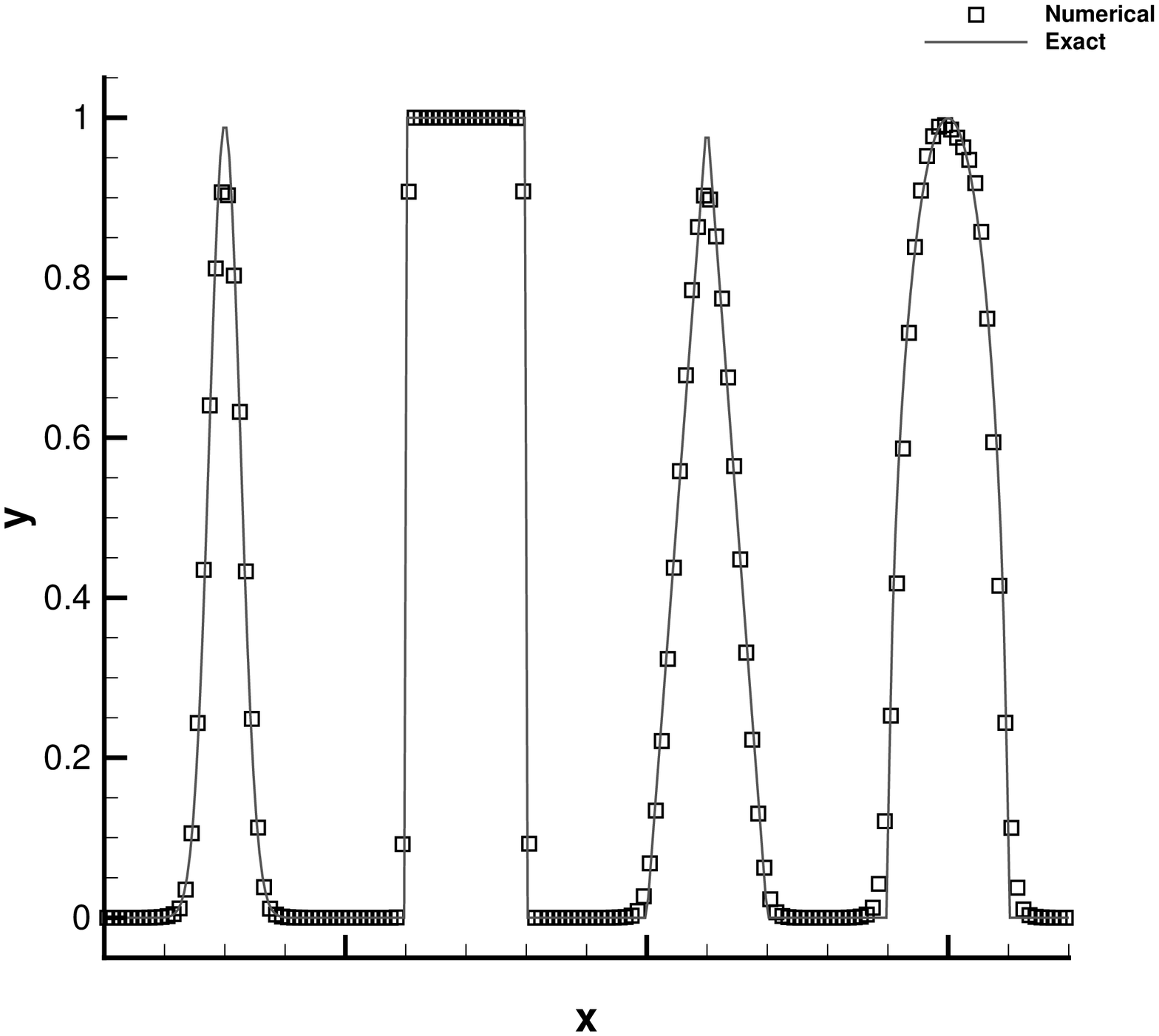}}
	\protect\caption{Same as Fig.\ref{fig:deng2wenoz}, but with $\beta=4.0$.\label{fig:deng2wenozBeta4}}	
\end{figure}

The numerical results calculated by the WENOZ-THINC-BVD scheme of different BVD algorithms have been shown in Fig.~\ref{fig:sunwenoz}\textendash\ref{fig:deng2wenoz}. All of them can solve discontinuities sharply by nearly 4 cells, which is a significant improvement of present schemes in comparison with those which only use high order polynomials in reconstructions. The result calculated by BVD(\rom{1}) is presented in Fig.~\ref{fig:sunwenoz}. Besides the discontinuous region, BVD(\rom{1}) also changes the solution around some critical points compared with the original WENOZ scheme. For algorithm BVD(\rom{2}) as shown in Fig.~\ref{fig:deng1wenoz}, the numerical result in the smooth regions look almost the same as the original one in Fig.~\ref{fig:wenoz}, while the numerical dissipation around the discontinuities is remarkably reduced.
 As one of algorithms devised for unstructured grids, algorithm BVD(\rom{3}) is capable of solving discontinuities sharply but pollutes smooth regions as shown in Fig.~\ref{fig:xiewenoz}, which may be caused by the assumption that neighbor cells are always smooth. Considering that the discontinuity cannot be resolved by only one cell, BVD(\rom{4}) is devised by evaluating TBV with the interpolation function over a group of neighboring cells. Shown in Fig.\ref{fig:deng2wenoz}, good results comparable to BVD(\rom{1}) and BVD(\rom{2}) can be obtained by BVD(\rom{4}) which is simple and can be directly implemented on unstructured grids. It is noteworthy that BVD(\rom{4}) can use even larger $\beta$ values. In Fig.~\ref{fig:deng2wenozBeta4}, we show the result with $\beta=4.0$ which generates sharper discontinuities resolved by only 2 cells.   
  
\section{Conclusion remarks \label{sec:conclusion}}

\section*{Acknowledgment}

This work was supported in part by JSPS KAKENHI Grant Numbers 15H03916, 15J09915 and 17K18838. 

\clearpage{}


\end{document}